\documentclass[intlimits,12pt,a4paper]{amsart}
\usepackage[all]{xy}
\usepackage{amssymb}

\newcommand{\re}{\operatorname{Re}}
\newcommand{\bt}{{\bf t}}
\newcommand{\bb}{{\bf b}}

\newcommand{\R}{{\bf R}}
\newcommand{\G}{{\bf G}}
\newcommand{\g}{{\bf g}}

\newcommand{\bH}{\bold{H}}
\newcommand{\bh}{\bold{h}}

\newcommand{\cO}{{\mathcal O}}

\newcommand{\A}{\mathbb{A}}

\newcommand{\Prod}{\prod\limits}

\newtheorem{theorem}{Theorem}[section]
\newtheorem{lemma}[theorem]{Lemma}

\numberwithin{equation}{section}
\begin{document}
\title[Conjecture of Jacquet and a local consequence]{Generalised form of a 
 conjecture of Jacquet and a local consequence}  
\author[D. Prasad \and  
R. Schulze-Pillot]{Dipendra Prasad \and  Rainer Schulze-Pillot} 
\thanks{This work was done during a stay of the second author
  at the Tata Institute of Fundamental Research, Mumbai. The authors thank
U.K. Anandavardhanan, Steve Kudla and Nimish Shah for helpful comments.} 
\maketitle
\begin{abstract}
{Following the work of Harris and Kudla we prove a general form
of a conjecture of Jacquet relating the non-vanishing of a certain
period integral to non-vanishing of the central critical value of a 
certain $L$-function. As a consequence, we deduce a theorem 
relating the existence of $GL_2(k)$-invariant linear forms on 
irreducible, admissible representations of $GL_2({\Bbb K})$ for 
a commutative semi-simple cubic algebra ${\Bbb K}$ over a non-archimedean
local field $k$ in terms of  local epsilon factors which was 
proved only in some cases by the first author in his earlier
work in \cite{dprasad_amj}. This has been achieved by globalising a 
locally distinguished representation to a globally distinguished 
representation, a result of independent interest which is proved by 
an application of the relative trace formula.}

\end{abstract}

\section{Introduction} 

Let $F$ be a number field  and $E$  a
semisimple cubic algebra over 
$F$, i.e., let $E$ 
be either a cubic field extension of $F$, or $E=F\oplus F'$ with a
quadratic field extension $F'$ of $F$, or $E=F\oplus F\oplus F$.

For a cuspidal automorphic representation $\Pi$ of
the group $GL_2(\A_E)$ with trivial central character restricted to 
$\A_F^{\times}$, where $\A_F$ is the ad\`ele ring of $F$, and
$\A_E= E \otimes_F \A_F$,
the period integral 
$$\int_{\A_F^{\times}GL_2(F)\backslash GL_2({\A_F})} f(g)d g,$$
has been much studied in the split case $E = F \oplus F \oplus F$, where
a famous question, asked by Jacquet, and now completely
settled by Harris and Kudla in \cite{hk_2} completing 
their earlier proof for special cases given in \cite{hk_1},
relates the non-vanishing of this period
integral to the non-vanishing of the central critical $L$-value
$L(\frac{1}{2}, \Pi_1\otimes \Pi_2 \otimes \Pi_3)$ of an 
automorphic representation $\Pi$ of $GL_2(\A_E) = GL_2(\A_F) {\times} 
  GL_2(\A_F) {\times} GL_2(\A_F)$ of the form $\Pi_1\otimes \Pi_2 \otimes \Pi_3$.

It is natural to study this question in the general case of 
any semisimple cubic algebra $E$  over a number field $F$. We recall that
for an automorphic representation $\Pi$ of $GL_2(\A_E)$, unramified outside a 
finite set $S$ of places of $E$, the triple product 
$L$-function comprising of local factors outside $S$,
to be denoted by $L_3^S(s,\Pi)$, was shown to have meromorphic
continuation and a certain functional equation by Piatetski-Shapiro
and Rallis in \cite{psr}, using a generalized version of Garrett's
integral representation. In \cite{ikeda2}, Ikeda proved 
(again using the integral
representation) that this result can be extended to the completed $L$-
function defined as a product of local factors over all places
(including the archimedean ones). The local factors occurring here are
defined in terms of local zeta integrals. 

On the other hand, Shahidi
identified  $L_3^S(s,\Pi)$ 
as a partial  $L$-function which arises out of 
the Langlands-Shahidi method, and
was able to get a complete $L$-function $L_3(s,\Pi)$
incorporating all the finite and infinite primes, having a functional
equation of the form 
$$ L_3(s, \Pi) = \epsilon_3(s,\Pi) L_3(1-s, \Pi),$$
where
 $$\epsilon_3( s, \Pi) 
= \prod_v
\epsilon_{3,v}( s, \Pi_v ),$$
and where $\Pi= \otimes \Pi_v$, $\Pi_v$ being a representation of 
$GL_2(E\otimes_F F_v)$, and $ \epsilon_{3,v}( s, \Pi_v )$ are
the local epsilon factors associated to $\Pi_v$ which are equal to 1
for almost all places $v$ of $F$.

\vspace{0.3cm}
If $D$ is a quaternion algebra over $F$ we denote by $D_E$ the
algebra $D\otimes_F E$ and by $\Pi^{D_E}$ the automorphic representation 
of $D^{\times}({\Bbb A}_E)$ corresponding to $\Pi$ under the correspondence of
Jacquet and Langlands (if it exists).

The aim of this paper is to prove the following general form of the
theorem of Harris and Kudla and to show that this  general form
can be applied  to the study of invariant linear forms, completing the results
of the first author in \cite{dprasad_amj}.

\begin{theorem}\label{jacquetconjecture_theorem}
Let $\Pi$ be a cuspidal irreducible automorphic representation of the adelic
group $GL_2({\Bbb A}_E)$ with central character that is trivial on 
${\Bbb A}_F^{\times}$.
Denote by $L_3(s,\Pi)$ the (completed) twisted tensor $L$-function of $\Pi$.
Then $L_3(\frac{1}{2},\Pi)\ne 0$ if and only if there exists a
quaternion algebra $D$ over $F$ such that $\Pi^{D_E}$ is defined and 
there exists a function 
%% $$F^D =\begin{cases}
%% f^D_1\otimes f^D_2\otimes f^D_3\\
%% f_1^D\otimes f_2^D\\
%% f^D
%% \end{cases} $$ 
$f^{D_E}$ in the space of $\Pi^{D_E}$ such that
\begin{equation}\label{conjecture_integral}
I(f^{D_E})=\int_{{\Bbb A}_F^{\times}D^{\times}(F)\backslash D^{\times}(\A_F)} f^{D_E}(y)d
y\ne 0.
\end{equation}
\end{theorem}

\vspace{2mm}

We show in this paper that most of the proof given in \cite{hk_2}
for the split case in fact goes through
for a general semisimple cubic algebra. It should be noted that
the proof in \cite{hk_2} uses the definition of the triple product
$L$-function through the integral representation from \cite{psr,
  ikeda2} and the fact that 
by the work of
Ramakrishnan  \cite{ramakrishnan_annals152}, 
%\marginpar{\em changed} 
this  $L$-function in the split case
 coincides with the definition of Shahidi as well
as with the one 
given in terms of representations of the
Weil-Deligne group; this has not yet been
established for general $E$ (in particular, 
Ikeda's paper \cite{ikeda2} excludes the case of a cubic
field extension $E/F$ in
Lemma 2. 2 and its
corollary).

We will see, however, that for questions of vanishing or nonvanishing 
at $1/2$
it makes no difference which definition of the $L$-function
$L_3(s,\Pi)$ we choose to
formulate our theorem. 

We end this introduction by mentioning that in the split case $E = F
\oplus F \oplus F$ more precise versions of Theorem
\ref{jacquetconjecture_theorem} have been investigated, 
starting with the work of Harris and Kudla \cite{hk_1} and continued by
Gross and Kudla \cite{gk}, Boecherer and Schulze-Pillot \cite{bs-triple}, 
and in the thesis of Watson \cite{WAT}. These more precise versions use the
connection between the central critical value of the $L$-function
$L_3(s,\Pi)$ with period integrals to prove rationality results 
and even explicit formulas for this central critical value under 
additional assumptions on $\Pi$; the explicit formulas obtained have
found applications to questions in quantum chaos in
\cite{WAT,boe-sa-sp}.  It would
be natural to extend these results in the context of our work to the
case of general cubic algebras over a number field.

% that not only the 
% nonvanishing properties of the period integral,
% $$\int_{\A_F^{\times}GL_2(F)\backslash GL_2({\A_F})} f(g)d g,$$
% but rationality questions, and even its precise value 
% (for a chosen {\it test vector f}, such as
% a new form) has been much studied in the split case 
% $E = F \oplus F \oplus F$, starting with the work of Harris and Kudla,
% Gross and Kudla, Boecherer-R. Schulze-Pillot, 
% and  in a very precise form in the thesis of Watson,
% which has found striking applications to questions in quantum chaos, as well
% as in the work of Bernstein-Reznikov [journal of Diff. Geometry, 2005]
% on subconvexity theorem for the triple product $L$-functions. It would
% be natural to extend these results in the context of our work to the
% case of general cubic algebras over a number field.

\section{Proof of Jacquet's conjecture}

Let $\G$ denote the algebraic group over $F$ whose group of $F$-points 
is the group  $$\{g \in GL_2(E)\mid \det(g) \in F^{\times}\}.$$ 

The group $\G(F)$ embeds  into the symplectic
similitude group $GSp_6(F)$ by viewing $GL_2(E)$ as the similitude
group of the standard symplectic form on $E^2$ which can be 
embedded in $GSp_6(F)$ by
choosing a $F$-basis of $E$ and taking the trace of the symplectic
form viewed as a symplectic form on $F^6$; the matrices from
$GL_2(E)$ with determinant in $F^{\times}$ yield then elements 
of the symplectic similitude group of this form
on $F^6$. This embedding carries over to the group of adelic points.
Let $\G_1$ denote the algebraic group over $F$ whose group of $F$-points 
is the group  $$\{g \in GL_2(E)\mid \det(g) = 1\}.$$ 

We consider the normalized Eisenstein series $E^{\times}(g,s,\Phi_s)$ on
$GSp_6(\A_F)$ defined
in \cite[(1.5)]{hk_2} with a factorizable section $\Phi_s$ and an
automorphic form $f \in \Pi$ 
%(which in the split cases 1) and 2) is assumed to be of the form  
%$f=f_1\otimes f_2\otimes f_3$ resp. $f=f_1\otimes f_2$)
and form the
global zeta integral
\begin{equation}
Z(s,f, \Phi_s)=\int_{{\Bbb A}_F^{\times}\G(F)\backslash\G(\A_F)}E^{\times}(\g,s,\Phi_s)
f(\g)d\g.
\end{equation}
We fix an additive character $\psi=\prod_v\psi_v$ of ${\A}_F$
and suppose that $f$ 
has a factorizable Whittaker function 
$W^\psi(g)=\prod_vW^\psi_v(g_v)$. 
%(which in the split cases 1) and 2) means
%that the $f_i$ have such Whittaker functions $W_i^\psi$).

By Theorem 2.1 of \cite{psr}, and the remark following it, 
the global zeta integral factors as a product  
$\prod_vZ_v(s,W_v^\psi,\Phi_{s,v})$
of local zeta integrals for sufficiently large $\re(s)$, 
where the local zeta integrals 
$Z_v(s,W_v^\psi,\Phi_{s,v})=\int\Phi_{s,v}(\delta_v g)W_v^\psi(g) dg$
are absolutely
convergent for large $\re(s)$ by \cite[Proposition 3.2]{psr}. This
factorization of the global zeta integral is obtained in \cite{psr}
(following Garrett's method) by the usual unfolding trick; the element
$\delta_v$ occurring in the last integral is a
representative of the only open orbit of
 $\G(F_v)$ acting on $P_v\backslash GSp_6(F_v)$, for $P_v$ the Siegel 
parabolic in $GSp_6(F_v)$.

The following Lemma is a  generalization of Corollary 2.3 of \cite{hk_2}:
\begin{lemma}\label{zetaintegral}
With notations as above one has $L_3(\frac{1}{2},\Pi) \ne 0$ if and only
if there is $f \in \Pi$ and a factorizable section $\Phi_s$ such that 
$Z(0,f,\Phi_s) \ne 0$.
\end{lemma}
\begin{proof} 
We prove the Lemma using the local factors given by the method of
integral representation and discuss the other definitions for the
local factors at bad places in the end.

We let $S$ be any finite set of places of $E$ containing all the archime\-dean
places such that the additive character $\Psi$ of $\A_E$ obtained from
the character $\psi$ of $\A_F$ by taking the trace map is unramified
outside $S$ and such that $\Pi$ is unramified outside $S$.

For a factorizable section $\Phi_s$ and $f \in \Pi$ with factorizable
 Whittaker function $W^\psi$ as above, we let $S_1=S_1(f,\Phi_s)\supseteq S$
be a finite set of places of $F$ such that 
for $v
\not\in S_1$, $f$ is  
$GL_2(\cO^E_v)$ fixed with  $W_v^\psi(e)=1$, where $\cO^E_v$ is the
maximal compact subring of $E\otimes F_v$, and such that the section
$\Phi_{s,v}$ is right 
$GSp_6(\cO_v)$-invariant with $\Phi_{s,v}(e)=1$.

By Theorem 3.1 of \cite{psr}  the local zeta integrals for the places
 $v \not\in S_1$ are then equal to the local factors of the $L$-function
$L_3(s,\Pi)$, with the argument shifted by $1/2$, so that we obtain
\begin{equation}
Z(s,f,\Phi_s)=L_3^{S_1}(s+\frac{1}{2},\Pi)\prod_{v \in
  S_1}Z_v(s,W_v^\psi,\Phi_{s,v}).
\end{equation}

This identity, at first valid for $\re(s)$ sufficiently large, extends
to the complex plane by meromorphic continuation.

As a first step, we prove the assertion of the Lemma with $L_3(s,\Pi)$
replaced by the partial $L$-function  $L_3^S(s,\Pi)$, where the
exponent $S$ on $L_3$ indicates the product of the local 
factors for the $v\not\in S$.

If  we have $L_3^{S}(\frac{1}{2},\Pi) \ne 0$ we can choose $f$
and the components $\Phi_{s,v}$ of a section $\Phi_s$ for the $v \not
\in S$ so that we can take $S_1=S$ above and have 
$$L_3^{S}(s+\frac{1}{2},\Pi)=\prod_{v \not\in
  S}Z_v(s,W_v^\psi,\Phi_{s,v}).$$
By 
\cite[Proposition 3.3]{psr}, \cite[p.227]{ikeda2} we can then choose the
  local data at the places in $S$ in such a way that the local zeta
  integrals $Z_v(s,W_v^\psi,\Phi_{s,v})$ for $v \in S$ are $\ne 0$ 
at $s=0$.
With this choice of data $f,\Phi_s$ we obtain then 
$$Z(0,f,\Phi_s)=\prod_{v \in
  S}Z_v(0,W_v^\psi,\Phi_{s,v}) \prod_{v \not\in
  S}Z_v(0,W_v^\psi,\Phi_{s,v})\ne
0,$$ since both factors are nonzero.

If $L_3^{S}(\frac{1}{2},\Pi) = 0$ then for any choice of
factorizable data we consider again the factorization
$$Z(s,f, \Phi_s)=\prod_vZ_v(s,W_v^\psi,\Phi_{s,v})$$ and choose a
finite set $S_1\supseteq S$ as above.

For the $v \in S_1\setminus S$
by \cite[App. 3 to Section 3]{psr},\cite[p. 227]{ikeda2}
the local zeta integral $Z_v(s-1/2,W_v^\psi,\Phi_{s,v})$ is a multiple of
the local $L$-factor  $L_{v}(s,\Pi)$ by a polynomial in
$q_v^{s},q_v^{-s}$, so that  
$$\prod_{v\not \in S}Z_v(s,W_v^\psi,\Phi_{s,v})$$ has a zero at $s=0$
as well. Moreover, by \cite[Lemma 2.1]{ikeda2} and using the result of 
\cite{kimshahidi_crparis} for the Satake parameters of
$\tilde{\Pi}_v$ at 
any place $v$, of an irreducible cuspidal
automorphic representation $\tilde{\Pi}$  of $GL_2$,  we see 
that the local zeta integrals $Z_v(s,W_v^\psi,\Phi_{s,v})$ 
have no pole in the region
$\re{s}>\frac{15}{34}-\frac{1}{2}$ (which contains the point $s=0$), 
for any finite $v$.

The zero of $\prod_{v\not \in
  S}Z_v(s,W_v^\psi,\Phi_{s,v})$ at $s=0$ 
therefore cannot be cancelled by a pole in the other factors and we
obtain $Z(s,f, \Phi_s)=0$ at $s=0$ for all choices of factorizable
  data as asserted.  

To come back to $L_3(\frac{1}{2},\Pi)$ itself, we notice that the local
factors of $L_3(s,\Pi)$ are nonzero by definition. Moreover, by Ikeda's
argument for the 
local zeta integrals used above, the local factors of $L_3(s, \Pi)$ 
have no poles at
$s=\frac{1}{2}$.  

We therefore see that $L_3(\frac{1}{2},\Pi) \ne 0$ is equivalent to
$L_3^{S}(\frac{1}{2},\Pi) \ne 0$, and the assertion of the Lemma
follows. 

We finally note that the argument given above for the transition between
$L_3(\frac{1}{2},\Pi)$ and
$L_3^{S}(\frac{1}{2},\Pi)$,
goes through for the two other methods mentioned in the introduction for defining
local $L$-factors at the places in $S$, so that the statement of
the Lemma is true for these other definitions of local factors as well
(even though we have at present no proof that the definitions do in
fact coincide).

Indeed what we need is that 
the local factors of $L_3(s, \Pi)$ 
have no poles at
$s=\frac{1}{2}$.  
This is clearly true for the local factors
given in terms of representations of the Weil-Deligne group  by
the Kim-Shahidi estimate from \cite{kimshahidi_crparis} used above. 
For the local factors given by
Shahidi in \cite{shahidi_annals_132} it is also true, we check this
for the case of a (local) cubic field extension: If $\Pi_v$ is
tempered, it is Proposition 7.2 of  \cite{shahidi_annals_132}. In
the (hypothetical) case that 
$\Pi_v$ is not tempered, it is a subquotient of a representation
induced from a character $(\mu_1 \vert\quad\vert^t,\mu_2
\vert\quad\vert^{-t})$ of the Borel subgroup. In that case, Shahidi's
construction in \cite[Section 7]{shahidi_annals_132} reduces the
definition of the local factor  by induction to local factors for
characters.
The latter are just the usual factors of Hecke $L$- functions with
pole free region $\{ s \mid \re(s)>3t\}$, and the estimate $t<5/34$
guarantees again that the local factor of $L_3$ is free of poles in
that region. 
\end{proof}
It remains to check that in the general case one can convert the
assertion of the last 
lemma about the global zeta integral into the assertion of the Theorem
in the same way as in \cite{hk_2}. 

We recall some notation. We let $D$ be a quaternion algebra over $F$
and $\alpha \in F^{\times}$ and consider $D$ equipped with  the
quadratic form $q_\alpha(x)=\alpha n(x)$ (where $n(x)$ denotes the
reduced norm of $x$) as a quadratic space $V=V(D,\alpha)$ over $F$. We
let $H=H_V=GO(V)$ be the similitude group of $V$ (with similitude factor
$\nu$)  and $H_1=O(V)\subseteq H$.

We denote by $GSp^{+,V}_6(\A_F)$ the adelic group
of symplectic similitudes (over $F$) having positive similitude factor
at all the (necessarily real)
infinite places $v$ of $F$ where $D\otimes F_v$ is a division algebra.
For $g \in GSp^{+,V}_6(\A_F)$ and a test function $\varphi=\varphi^V 
\in S(V(\A_F)^3)$ the usual theta 
integral (defined by regularization if $D=M_2(F)$), is defined as
\begin{equation}\label{thetaintegral}
I(g,\varphi)=\int_{H_1(F)\backslash H_1(\A_F)}
\theta(h_1h,g;\varphi)dh_1,
\end{equation}
where $h \in H(\A_F)$ with $\nu(h)=\nu(g)$ as in
\cite[(3.9)]{hk_2}. Moreover, the function $I(\cdot,\varphi)$ on
$GSp^{+,V}_6(\A_F)$ has a unique left $GSp_6(F)$-invariant extension to
all of $GSp_6(\A_F)$, which we will also denote by  $g \mapsto
I(g,\varphi)$.

By \cite[Corollary 5.1]{hk_2} the assertion of Lemma
\ref{zetaintegral}
gives then,  as a  consequence of the
Siegel-Weil theorem, that one has
$L_3(\frac{1}{2},\Pi)=0$ if and only if for all choices of
isometry class of $V=V(D,\alpha)$ as $D$ and $\alpha \in F^*$ vary, test function $\varphi^V \in
S(V(\A_F)^3)$  and $f \in \Pi$ one has
\begin{equation}\label{thetacriterion}
\int_{\A_F^{\times}\G(F)\backslash\G(\A_F)}I(\g,\varphi^V)f(\g)d\g=0;
\end{equation}
in this step the structure of the cubic algebra $E$
plays no role at all.

For each choice of $V$, let  $\bH(F)=\bH_V(F) = \{h \in 
H(E)\mid \nu(h) \in F^{\times}\}$.
%$$\bH(F)=\bH_V(F)=\begin{cases}
%\{(h_1,h_2,h_3) \in H^3\mid \nu(h_1)=\nu(h_2)=\nu(h_3)\}\\
%\{(h_1,h_2) \in H {\times} H_{E'} \mid \nu(h_1)=\nu(h_2)\}\\
%\{h \in H_E\mid \nu(h) \in F\}
%\end{cases},$$
%where the distinction of cases as well as $E',E$ are as always in this
%section and where $H_{E'}$ resp. $H_E$ is the similitude group of the
%quadratic space $V\otimes E'$ resp. $V\otimes E$. 
We write $\bH_1$ for the subgroup of $\bH$ where the similitude
groups above are replaced by the respective orthogonal groups.

We write $\tilde{H}=(D^{\times} {\times}
D^{\times})\rtimes \langle \bt \rangle $, where the involution $\bt$ acts
on $D^{\times} {\times} 
D^{\times}$ by $(b_1,b_2) \mapsto (\bar b_2^{-1},\bar{b_1}^{-1})$,
and $\tilde{H^0}=D^{\times} {\times}
D^{\times}$.
Analogously we have  $\tilde{\bH}, \tilde{\bH}^0$ as covering groups of
$\bH=GO(V\otimes E)$, and
$\bH^0 = GSO(V\otimes E)$.

The groups $Sp_6, H_1$ form  a dual reductive pair in $Sp_{24}$,
and the groups $\G_1, \bH_1$ form another dual reductive pair in this
big symplectic group.
For the  latter pair,  we obtain first in the usual way an
embedding of $SL_2(E){\times} O(V\otimes E)$ into the isometry group
of the 
standard $8$-dimensional symplectic space over
$E$  which is then transferred to an embedding of  $\G_1 {\times}  \bH_1$
into $Sp_{24}$ by taking the trace from $E$ to $F$.

The proof in \cite{hk_2} proceeds then (for $E=F^3$)
as follows: For a cuspidal automorphic form 
$f$ on $\G(\A_F)$ one considers the lifted form on $\bH(\A_F)$ given by
\begin{equation}
I(\bh, \varphi;f) 
=\int_{\G_1(F)\backslash \G_1(\A_F)}\theta (\bh,
\g_1\g;\varphi)f(\g_1\g)d\g_1, 
\end{equation} where $\g \in \G(\A_F)$ is arbitrary with
$\det(\g)=\nu(\bh)$ and where $\varphi$ is any Schwartz-Bruhat
function on the space $V(\A_F)^3$. 

One considers then the seesaw pair 
\begin{equation}
\xymatrix{{I(\cdot, \varphi)}& GSp_6 
\ar@{-}[d]\ar@{-}[dr]&\bH\ar@{-}[d]& I(\cdot, \varphi, f)\\
  f & \G \ar@{-}[ur]& H=GO(V)& 1 \quad \quad\quad,} 
\end{equation}
where $\bH=\bH_V$ and $H=H_V$ depends on
the choice of the 
quadratic space $V=V(D,\alpha)$). 

The seesaw identity given in \cite{kudla_seesaw}
yields in this situation
\begin{equation}\label{seesaw}
\int_{\A_F^{\times}\G(F)\backslash \G(\A_F)}I(\g,\varphi^V)f(\g)d\g=
\int_{\A_F^{\times} H(F)\backslash H(\A_F)}I(h,\varphi^V,f)dh
\end{equation}
for all the possible choices of $V=V(D,\alpha)$. 

From this we see that $L_3(\frac{1}{2},\Pi)=0$ if and only if for 
all choices of $V=V(D,\alpha)$, all
$\varphi \in S(V(\A_F)^3)$, and all $f \in \Pi$,  
one has
\begin{equation}\label{integralcrit}
\int_{\A_F^{\times} H(F)\backslash H(\A_F)} I(h,\varphi^V,f)dh=0.
\end{equation}

Finally, fixing one choice of $V=V(D,\alpha)$ (and
omitting  $V$ from
the notation) it is shown that 
\begin{equation}\label{uniqueness} 
\int_{\A_F^{\times}H(F)\backslash H(\A_F)}I(h,\varphi,f)dh=
\int_{(\A_F^{\times} \times \A_F^{\times})\tilde{H}^0(F)\backslash
  \tilde{H}^0(\A_F)}\tilde{I}(h,\varphi,f)dh, 
\end{equation}
where 
$\tilde{I}(\cdot,\varphi,f)$ denotes the
pullback of $I(\cdot,\varphi,f)$ to
$\tilde{\bH}(\A_F)$. Rewriting
\begin{equation}\label{restriction}
\tilde{I}((\bb_1,\bb_2),\varphi,f)=
\sum_rI^{1,r}(\bb_1,\varphi,f)I^{2,r}(\bb_2,\varphi,f)   
\end{equation}
for $\bb_i \in D^{\times}(\A_E)$ and for suitable 
functions $I^{i,r}(\bb_1,\varphi,f) \in \Pi^D$ one arrives at the
assertion of the Theorem. (We recall that an automorphic function
on a product of two groups is a finite sum of a product of automorphic
functions on the two groups.)

The key
tools in \cite{hk_2} in these steps are the analysis of the theta
correspondence for 
similitude groups for the dual pair $(\G=(GL_2(F)^3)_0, \bH=(GO(V)^3)_0)$ 
where the subscript 
0 refers to the subgroup with equal similitude factors) and the use of the
seesaw identity given above. 
Since the necessary facts about the theta correspondence for
similitude groups  were proved in \cite{hk_2} over an arbitrary number field,
this analysis still applies if one replaces the pair
$(GL_2(F)^3,D^{\times}(F)^3)$  by  $(GL_2(E),D^{\times}(E))$ or
$(GL_2(F){\times} GL_2(F'),D^{\times}(F){\times} D^{\times}(F'))$.
It is then not difficult to check step by step that the proof carries
over to our more general situation. 
For the sake of completeness we conclude this section by reconsidering
two of the steps given above.

%%%%%paragraph changed June 10%%%%%%%%%%
The first point that we wish to elaborate upon is about the case
$D(F)=M_2(F)$, 
where as remarked earlier the theta 
integral giving the function $I(g,\varphi)$ on $GSp_6(\A_F)$ has to be
defined by regularization. If $F$ has at least one real place $v_0$
this can
be done as described in \cite[p. 620/621]{hk_1}.

One considers  then the differential operator $z$ from the center of the
universal enveloping algebra of $GSp_6(F_{v_0})$ (with $F_{v_0}=\R$)
whose image under the Harish-Chandra isomorphism is the symmetric
polynomial
$(x_1^2x_2^2+x_1^2x_3^2+x_2^2x_3^2-5(x_1^2+x_2^2+x_3^2)+21)/12$. 
This  differential  operator has the property that 
it maps to zero
in the universal enveloping  algebra of $SO(1,1)$ under the
homomorphism of universal enveloping
algebras induced by the local theta correspondence, from which one can
deduce that 
for any Schwartz function 
$\phi \in {\mathcal S}(M_{2,6}({\Bbb R}))$ 
with Fourier transform $\hat{\phi}\in {\mathcal S}(M_{6,2}({\Bbb R}))$, 
$\omega(z)\hat{\phi}$  vanishes on the set of matrices of rank $\leq 1$. As a 
consequence, $\theta(g,h,\omega(z)\phi)$ is of rapid decay along 
$O(V)(F)\backslash O(V)({\Bbb A}_F)$, and hence it makes sense to
define
\begin{equation}
I(g,\varphi):=\int_{ H_1(F)\backslash
  H_1(\A_F)}\theta(g,h_1h,\omega(z)\varphi)dh_1
\end{equation} 
where $h \in H(\A_F)$ with $\nu(h)=\nu(g)$ as above.

Using this modified $I(g,\varphi)$ and writing $V=V_0$ in this case
one obtains  
\begin{equation}\begin{split}
\int_{\A_F^{\times}\G(F)\backslash\G(\A_F)}&I(\g,\varphi^{V_0})f(\g)d\g\\
&=\int_{\A_F^{\times}\G(F)\backslash\G(\A_F)}\Bigl(\int_{ H_1(F)\backslash
  H_1(\A_F)}\theta(\g,h_1h,\omega(z_H)\varphi)dh_1\Bigr)f(\g)d\g\\
&=
\int_{\A_F^{\times}H(F)\backslash H(\A_F)}I(h,\omega(z_H)\varphi^{V_0},f)dh \\
&=\int_{\A_F^{\times}H(F)\backslash
  H(\A_F)}r(z_H)I(h,\varphi^{V_0},f)dh\\
&= c\int_{
\A_F^{\times}H(F)\backslash  H(\A_F)
}I(h,\varphi^{V_0},f)dh
\end{split}
\end{equation}

Here, in the first equality we take an element $z_H$ of the center of
the universal  enveloping
algebra of the Lie algebra $\frak{h}$ of $H(\R)$ with
$\omega(z)=\omega(z_H)$ on the space $S(V_0(\A_F)^3)$; such an element
exists and can be computed explicitly as described in 
\cite{przebinda},\cite[Proposition 5.1.1]{kudla_rallis_regularized}. 
The second equality is 
 the usual seesaw identity with the differentiated theta kernel
 $\theta(g,h,\omega(z_H)\varphi)$ which is of rapid
decay.
The third equality, in which we write $r(z_H)$ for the action of $z_H$ on
functions on $H(\A_F)$, is a consequence of the fact that 
the map $\varphi \mapsto I(h,\varphi,f)$ from  $S(V_0(\A_F)^3)$ to
the space  
of  functions on 
$ \A_F^{\times}H(F)\backslash  H(\A_F)$ is $H(\A_F)$-invariant.  

Since the theta lift $I(h,\varphi^{V_0},f)$ of the cusp form $f$ to
$H(\A_F)$ is cuspidal we obtain the last equality (with $c$ being equal to
the constant term of $r(D_H)$) from the well known fact that 
$$\int_{ \A_F^{\times}H(F)\backslash H({\Bbb A}_F)} \phi_1(h)X \phi_0(h) dh =
-\int_{ \A_F^{\times}H(F)\backslash H({\Bbb A}_F)} X\phi_1(h) \phi_0(h) dh$$
for $X \in \frak{h}$ and a pair consisting of a  cusp form $\phi_0$ and an
arbitrary  automorphic form $\phi_1$ on
$H(\A_F)$ with trivial central characters.

An explicit computation gives $c=1$, and one obtains an identity of the
usual shape 
$$\int_{\A_F^{\times}\G(F)\backslash\G(\A_F)}I(\g,\varphi^{V_0})f(\g)d\g=\int_{
\A_F^{\times}H(F)\backslash  H(\A_F)
}I(h,\varphi^{V_0},f)dh.$$

Thus although the seesaw identity is invoked only for certain elements of
 the space $S(V_0(\A_F)^3)$ of the form $\omega(z) \varphi$, these 
elements  lift cusp forms on $\G$ surjectively onto cusp forms on ${\bf H}$. 

If all archimedean places of $F$ are complex one can replace the
argument above by a similar argument using a nonarchimedean place at
which all representations are unramified principal series and
replacing the differential operator used above by a suitable Hecke
operator as in \cite{tan}. Explicitly, one uses the Hecke operator
$z$ whose image under the Satake isomorphism is the polynomial
$$x_1x_2+x_1x_3+x_2x_3-(y_2+y_1)(x_1+x_2+x_3)+(y_4+y_3+y_2+y_1+4)$$
where $x_j=q^{s_j}+q^{-s_j},y_j=q^j+q^{-j}$ and where $q$ is the order
of the residue field. This $z$ maps to zero
in the Hecke algebra of $SO(1,1)$ under the homomorphism of Hecke
algebras induced by the local theta correspondence and can therefore
be used for 
regularizing the theta integral (similarly as done above at the real
place) as in \cite{tan}; its image in the
Hecke algebra of $SO(2,2)$ acts by the nonzero scalar $y_4-y_3-2y_2+y_1$ on the
constant function $1$, and the argument from above works in the same way. 

\vspace{0.2cm}
%%%%%end of paragraph changed June 10%%%%%%%%%%
Our second remark concerns the proof of (\ref{uniqueness}) in
\cite{hk_1,hk_2} relating the period integral on $GO(4)$ to one on  $GSO(4)$
 which we elaborate upon in the next section.

\section{Automorphic forms on disconnected groups}

The period integral which naturally occurs from the
consideration of seesaw pairs in the last section 
is on the group $GO(4)$, whereas the final theorem
is about period integral on $GSO(4)$. (By 2.10, it is
trivial to see that the integral on $GSO(4)$ is nonzero if and only if 
the integral on the corresponding quaternion 
algebra is nonzero.) 
Harris and
Kudla give a proof of this in [6]. This is a  subtle point. 
We take the occasion to make some general comments about automorphic forms on disconnected groups. 
We begin with some generality.

Let $G$ be a  group containing a subgroup $G_0$ of index 2, and
let $\sigma$ be an involution in $G$ such that $G = G_0 
\rtimes \langle \sigma \rangle$. 

An irreducible representation of $G_0$ which is invariant under $\sigma$
can be extended in exactly two distinct 
ways to an irreducible representation of $G$.
In general, there is no way of distinguishing between these two extensions. However,
there are many situations in which one can
distinguish between these two extensions. One context in which one has 
learnt to distinguish between the two extensions is in the problem
of base change for $GL_n$, say over a local field. Here given a 
quadratic extension $K/k$ of local fields with $\sigma$ 
the nontrivial element of the Galois group, for an irreducible admissible 
representation
$\pi$ of $GL_n(K)$ which is Galois invariant, there is a natural extension
$\tilde{\pi}$ of $\pi$ to $GL_n(K) \rtimes \langle \sigma \rangle$ such that
the character identity,
$$\Theta_{\tilde{\pi}}(g,\sigma) = \Theta_\tau({\rm Nm} g),$$
holds for a representation $\tau$ of $GL_n(k)$, $g$ any $\sigma$-regular
element of $GL_n(K)$ with ${\rm Nm} g$ its norm.
This natural extension is achieved via uniqueness of Whittaker, or
degenerate Whittaker model through a generality that we discuss below.

In general, suppose $G_0$ has a subgroup $H$ together with a character
$\chi: H \rightarrow {\Bbb C}^{\times}$ such that the inner-conjugaction action by 
$\sigma$ takes $H$ to 
$H$, and 
$\chi$ to $\chi$. Assuming that $\chi$ appears in $\pi$ with multiplicity 1,
we define a {\it preferred} extension $\tilde{\pi}$ 
of $\pi$ to be the representation of 
$G$ for which the character
$\tilde{\chi}: H \rtimes \langle \sigma \rangle \rightarrow {\Bbb C}^{\times} $
extended from $\chi$ trivially across $\langle \sigma \rangle$ 
appears (as a quotient) in
$\tilde{\pi}$.  Of course this notion of preferred extension depends on the $G$-conjugacy
class of the triple $(\sigma, H,\chi)$. It happens in many situations that there is more
than one natural choice of $G$-conjugacy class of the triple $(\sigma,H,\chi)$ in which case 
it may be of interest to compare the various preferred extensions.

In our context, there exists a canonical extension of a
representation of the group $GSO(4)$ over a local field 
to one of $GO(4)$ which we discuss now.

Let $(q,V)$ be a quadratic space of dimension 4, and discriminant 1. It is easy to
see that $GO(V)$ operates transitively on the set of vectors $V_{q \not = 0}=\{v \in V| 
q(v) \not = 0\}$. Fix $v_0 \in V$ such that $q(v_0) \not = 0$, and let $W$ be the orthogonal
complement of $v_0$, and $\sigma$ the reflection around $W$. Let $H=O(W)$ be the stabiliser 
in $GO(V)$ of the vector $v_0$. It can be seen that an irreducible representation $\pi$ 
of $GSO(V)$ which is invariant under the action of $\sigma$ contains a unique, up to 
scaling, linear form $\ell: \pi \rightarrow {\mathbb C}$ fixed under $H_0 = SO(W)$. The
preferred extension of $\pi$ to $GO(V)$ is the one for which $H = O(W)$ operates
trivially on $\ell$. 

It can be seen that if $(q,V)$ is the direct sum of 2 Hyperbolic spaces, then 
for an irreducible unramified representation $\pi$ of $GSO(V)$, the preferred extension
as defined in the previous paragraph is the unique spherical representation of $GO(V)$
containing $\pi$.

Similarly  there exists a canonical extension of a $\sigma$-invariant 
automorphic representation of the group $GSO_4({\Bbb A}_F)$  
to one of $GO_4({\Bbb A}_F)$ which we describe now, in greater generality.

Suppose $G$ is an algebraic group over a number field $F$, containing an involution $\sigma$, 
and a subgroup $G_0$ such that $ G = G_0 \rtimes <\sigma>$. Let $\pi = \otimes \pi_v$ 
be an automorphic representation of $G_0$
realised on a space of functions on $G_0(F)\backslash G_0({\Bbb A}_F)$ which
is invariant under $\sigma$. Since the factorisation $\pi = \otimes \pi_v$ 
is unique up to isomorphism, each of the representations $\pi_v$ 
is invariant under $\sigma_v$, i.e., $\sigma_v(\pi_v) \cong \pi_v$; this cannot be used 
to define an action of $G(F_v)$ on $\pi_v$ as the action of $\sigma_v$ is well-defined only up to a 
sign.

Recall that evaluating at the identity element gives rise to a linear form
$e_\pi: \pi \rightarrow {\mathbb C}$ which is $G_0(F)$-invariant. This linear
form is clearly $\sigma$-invariant for the natural action of $\sigma$ on $\pi$.
To extend the automorphic representation $\pi = \otimes_v \pi_v$ of $G_0({\Bbb A}_F)$
to $G({\Bbb A}_F)$, first we extend the representation $\pi_v$ of $G_0(F_v)$
to a representation $\tilde{\pi}_v$ of $G(F_v)$, keeping it spherical at almost all places, 
and such that the 
action of $\sigma$ on $\pi = \otimes_v \pi_v$ is the same as that of $\otimes_v \sigma_v
$ on it, but with no other
constraints. (We will call the condition $\sigma = \otimes_v \sigma_v$, 
the {\it coherence} condition.) We keep the linear form $e_{\tilde{\pi}}: \tilde{\pi}=\otimes_v 
\tilde{\pi}_v \rightarrow {\Bbb C}$ 
   same as $e_\pi$ on $\pi = \tilde{\pi}$. Since $e_\pi$ is $\sigma$ as well as $G_0(F)$-invariant,
it is $G(F)$-invariant, and thus we have constructed an automorphic representation
of $G({\Bbb A}_F)$. 

In particular, if an automorphic representation $\pi = \otimes_v \pi_v $ of $G_0({\Bbb A}_F)$ is 
$\sigma$-invariant, and the 
representations $\pi_v$ of $G_0(F_v)$
have a canonical extension  to representations $\tilde{\pi}_v$ of 
$G(F_v)$, there is a canonical
extension of automorphic forms from $G_0({\Bbb A}_F)$ to $G({\Bbb A}_F)$
assuming that the coherence condition $\sigma = \otimes_v \sigma_v$ is satisfied.
In the context of $GSO(4)$, this
canonical extension is the one that appears as theta lift from $GL(2)$,
cf. \cite {hps}, the coherence condition being automatic here.

The following lemma relates the non-vanishing of period integrals on $GO_4({\Bbb A}_F)$
to that on $GSO_4({\Bbb A}_F)$.

\begin{lemma}
Let $E$ be a cubic semi-simple algebra over a number field $F$, and $\Pi$ an automorphic
representation on $GO_4({\Bbb A}_E)$ obtained by theta lifting from $GL_2({\Bbb A}_E)$. 
Then a nonzero $GSO_4({\Bbb A}_F)$-invariant linear form on $\Pi$ is $GO_4({\Bbb A}_F)$-invariant.
\end{lemma}

\noindent{\bf Proof :} 
Because of the multiplicity 1 theorems in 
\cite{p-trilinear}, and \cite{dprasad_amj}, it suffices to prove the corresponding local statement. 
Let ${\Bbb K}$ be a cubic semi-simple algebra over a local field $k$, and
 $\pi$ an irreducible representation of $GO_4({\Bbb K})$ obtained by theta lifting from $GL_2({\Bbb K})$. 
Then we prove that any $GSO_4(k)$-invariant linear form on $\pi$ is $GO_4(k)$-invariant. 

Identifying $GSO_4$ to the quotient of $D^* \times D^*$ 
by the scalar matrices $(t,t^{-1})$, where $D$ is a quaternion algebra over $k$,
 $GO_4$ becomes the corresponding quotient of 
$(D^* \times D^*) \rtimes {\Bbb Z}/2$ where ${\Bbb Z}/2$ operates by switching the 2 factors. 

Under the identification of $GSO_4$ with $ D^* \times D^*$ divided by a central subgroup, the
representation $\pi$ which is invariant under the switching of the factors, can be written as
$\pi_1 \otimes \pi_1$ for an irreducible representation $\pi_1$ of $D^*({\Bbb K})$. 
Any $GSO_4(k)$-invariant
linear form on $\pi_1 \otimes \pi_1$ is up to a scalar, of 
the form $\ell \otimes \ell: \pi_1 \otimes  \pi_1 \rightarrow {\Bbb C}$ 
for a $D^*(k)$-invariant linear form $\ell: \pi_1 \rightarrow {\Bbb C}$. 
Clearly, the linear form $\ell \otimes \ell: \pi_1 \otimes  \pi_1 \rightarrow {\Bbb C}$ 
is invariant under  the switching of the factors, and that completes the proof of the lemma.

\section{Globalisation of locally distinguished representations}

In this section we prove a globalisation theorem which is  the main 
technical tool for establishing  local theorems from the 
corresponding global theorems.

Let $k$ be a local field, $G$ a reductive algebraic 
group over $k$, and $H$ a closed 
algebraic subgroup of $G$ defined over $k$. We abuse notation and 
use $G,H$ to also denote the corresponding group of $k$-rational points.  
For a character $\chi: H \rightarrow {\Bbb C}^{\times}$,
we define a 
representation $\pi$ of $G$ to be $\chi$-distinguished by $H$ if 
there exists
a nonzero linear form $\ell: \pi \rightarrow {\Bbb C}$ on which $H$
operates via the character $\chi$. 

Let $\underline{G},\underline{H}$ now be algebraic groups 
defined over a number field $F$, and 
let $\chi$ be a one dimensional
automorphic representation of $ \underline{H}({\Bbb A}_F)$. 
Let $\underline{Z}$ be the identity component of the center  of the 
algebraic group $\underline{G}$, and which we assume without loss of generality in the
rest of this section to be contained in $\underline{H}$. We assume further that 
$\underline{H}/\underline{Z}$ has no $F$-rational character.
We abuse notation to denote the restriction of $\chi$ to 
$\underline{Z}(F)\backslash \underline{Z}({\Bbb A}_F)$ also by $\chi$.
An automorphic representation $\Pi$ of $\underline{G}({\Bbb A}_F)$ 
is said to be globally $\chi$-distinguished by $\underline{H}$, if the period integral 
$$\int_{\underline{H}(F)\underline{Z}({\Bbb A}_F)\backslash \underline{H}({\Bbb A}_F)} f(h)\chi^{-1}(h) dh, $$
is nonzero for some $f \in \Pi$. Observe that for the function 
$ f(h)\chi^{-1}(h)$ on $\underline{H}({\Bbb A}_F)$ 
to be $\underline{Z}({\Bbb A}_F)$-invariant, we must have the central character of $\Pi$
restricted to  $\underline{Z}({\Bbb A}_F)$ be the same as $\chi$ restricted to 
$\underline{Z}({\Bbb A}_F)$. We note that by the theorem of Borel and Harish-Chandra,
under our assumption that $\underline{H}/\underline{Z}$ has no $F$-rational characters, 
$\underline{H}(F) \underline{Z}({\Bbb A}_F)\backslash 
\underline{H} ({\Bbb A}_F)$ has finite
volume. (Here $\underline{H}$ need not be reductive.)  Therefore the period integral 
makes sense, for instance, for $f$ a cusp form,
 as cusp forms are known to be bounded 
on $\underline{G}({\Bbb A}_F)$. 

Here is the main theorem of this section.

\begin{theorem}
Let $\underline{H}$ be a closed subgroup of a reductive group $\underline{G}$, both defined over 
a number field $F$. Let $\underline{Z}$ be the identity component 
of the center  of $\underline{G}$, which we assume is contained in $\underline{H}$ such that ${\underline{H}/\underline{Z}}$ 
has no $F$-rational character. 
Let $\chi = \prod_v \chi_v$ be a one dimensional
automorphic representation of $ \underline{H}({\Bbb A}_F)$.  
Suppose that $S$ is a finite set of non-Archimedean places of $F$, and
$\pi_v$ a supercuspidal representation of $\underline{G}(F_v)$ for all $v \in S$, 
which is $\chi_v$-distinguished
by $\underline{H}(F_v)$, i.e., ${\rm Hom}_{\underline{H}(F_v)}(\pi_v, \chi_v) \not = 0$ for all 
$v \in S$. Let $T$ be a finite set of places containing $S$ and all the 
infinite places,
such that $\underline{G}$ is  quasi-split at places outside $T$, 
and $\chi_v$ is unramified outside $T$, i.e., if $\underline{G}(\cO_v)$ is a 
hyperspecial maximal compact subgroup of $\underline{G}(F_v)$, then $\chi_v$ is trivial 
on $\underline{H}(F_v) \cap \underline{G}(\cO_v)$. 
Then there
exists a global automorphic form $\Pi = \otimes \Pi_v$ of $\underline{G}({\Bbb A}_F)$, necessarily
cuspidal, such that $\Pi_v = \pi_v$ for $v \in S$, and $\Pi_v$ is unramified
at all finite places of $F$ outside $T$, and an $f \in \Pi$
such that $$\int_{\underline{H}(F)\underline{Z}({\Bbb A}_F)\backslash \underline{H}({\Bbb A}_F)} f(h)\chi(h)^{-1} dh \not = 0.$$
\end{theorem} 

Before beginning the proof of this theorem, we gather together a few
results which go into the proof. The first result we need is 
a lemma about when certain global points in an adelic space are discrete.

\begin{lemma} 
Let $\underline{H}$ be a closed algebraic 
subgroup of a reductive group $\underline{G}$, both defined over 
a number field $F$. 
Assume that there exists a finite dimensional algebraic  representation $V$ 
of $\underline{G}$ defined over $F$, and a vector
$v \in V$ whose stabiliser in $\underline{G}$ is $\underline{H}$. 
Then the orbit of $\underline{G}(F)$ passing through the identity
element of $\underline{G}({\Bbb A}_F)/ \underline{H}({\Bbb A}_F)$ is a discrete subset of 
$\underline{G}({\Bbb A}_F)/ \underline{H}({\Bbb A}_F)$.
\end{lemma}
\noindent{\bf Proof :} The hypothesis of the lemma identifies
  $\underline{G}({\Bbb A}_F)/ \underline{H}({\Bbb A}_F)$ to a subset
of $V({\Bbb A}_F)$ containing the vector $v \in V(F)$. Since the $\underline{G}(F)$ orbit
of $v \in V(F)$ is contained in $V(F)$ which is a discrete subset of   $V({\Bbb A}_F)$,
the lemma follows. (We remark that if $\phi: X' \rightarrow X$ is a 
continuous injective map of topological spaces, $D$ a subset of $X'$ for which
$\phi(D)$ is discrete in $X$, then $D$ is discrete in $X'$. We use this
 remark for $X' = \underline{G}({\Bbb A}_F)/ \underline{H}({\Bbb A}_F)$, and $X =V({\Bbb A}_F)$.)

\vspace{2mm}

For the following well-known lemma about algebraic groups 
see for instance Propositions 7.8 and  7.7 of \cite{Borel}.

\begin{lemma}

{\bf 1.} Let $H$ be a closed subgroup of an algebraic  group $G$, 
both defined over 
a field $k$. Assume that $H$ has no nontrivial characters defined over $k$. 
Then there exists a finite dimensional algebraic  representation $V$ 
of $G$ defined over $k$, and a vector
$v \in V$ whose stabiliser in $G$ is $H$. In particular, if 
$Z$, the identity component 
of the center  of $G$,  is contained in $H$, and if $H/Z$ has no $k$-rational
characters,  then there exists a finite dimensional algebraic  
representation $V$  of $G$ defined over $k$, and a vector
$v \in V$ whose stabiliser in $G$ is $H$.

{\bf 2.} Let $H$ be a reductive subgroup of a connected algebraic group $G$, 
both defined over 
a field $k$. Then there exists a finite dimensional algebraic  
representation $V$ 
of $G$ defined over $k$, and a vector
$v \in V$ whose stabiliser in $G$ is $H$. 

\end{lemma}

For a continuous function $f$ on $G$ belonging to 
$C_c^\infty(Z\backslash G, \chi)$ (these are functions $f$ on $G$ with
$f(zg) = \chi(z)f(g)$ for $z\in Z$, $g \in G$, and compactly supported
modulo $Z$), define
the `orbital integral' along $H$ to be the function on $G$ given by
$$g \mapsto \Phi(f, g):= \int_{Z \backslash H}f(h g)
\chi(h)^{-1} dh.$$
Clearly $\Phi(f, g)$ is in the set  $C_c^\infty(H\backslash G, \chi)$
of  continuous functions
$\phi$ on $G$ which have the property that $\phi(hg) = \chi(h)\phi(g)$ 
for $h\in H$, $g \in G$ and are compactly supported
modulo $H$. 
It is easy to see that the `orbital integral' map is a surjection from 
$C_c^\infty(Z\backslash G, \chi)$  to $C_c^\infty(H\backslash G, \chi)$.

The following  general lemma is a consequence of Schur orthogonality relations,
cf. the Corollary in the appendix to \cite{hm}

\begin{lemma} Let $k$ be a non-Archimedean local field, $G$ a reductive 
algebraic group over $k$, 
and $H$ a closed
subgroup of $G$.  
Let $Z$ be the identity component
of the center  of $G$, which we assume is contained in $H$. Let 
$\chi: H \rightarrow {\Bbb C}^{\times}$ be a character of $H$ as well as its
restriction to $Z$. Let 
$(\pi,V)$ be a supercuspidal  representation of $G$, and
$(\pi', V')$ its contragredient. 
Suppose $\ell: \pi \rightarrow {\Bbb C}$ be a linear form on $\pi$ on which 
$H$ acts by $\chi$, and $v \in \pi $ is chosen so that $\ell(v) \not = 0$.
For  $v' \in V'$, let $f(g) = v'(gv)$ be a matrix coefficient of $\pi$. 
Then $$g \rightarrow \int_{Z \backslash H}f(hg) \chi(h)^{-1}dh$$
is a nonzero smooth function on $G$ with compact support in $H \backslash G$, i.e.,
belongs to $C_c^\infty(H\backslash G, \chi)$.
\end{lemma}

Finally, here is the statement of the relative trace formula that we will use;
see \cite{hm}, Theorem 2.

\begin{theorem}
Assume that $f = \otimes f_v \in C_c^\infty(\underline{Z}({\Bbb A}_F)\backslash 
\underline{G}({\Bbb A}_F), \chi)$ is such that at some finite 
place $v$ of $F$,
$f_v$ is the matrix coefficient of a supercuspidal representation $\pi_v$ 
of $\underline{G}(F_v)$. Then,
$$\sum_{\pi} \sum_{\phi \in {\Bbb B}_\pi} P(R(f) \phi) \overline{\phi(1)} = 
\sum_{\gamma \in \underline{H}(F) \backslash \underline{G}(F)} \Prod_v \Phi_v(f_v,\gamma),$$
where $\pi$ ranges over $\chi$-distinguished  cuspidal
representations of $\underline{G}({\Bbb A}_F)$ which have $\pi_v$ as the local 
component at $v$, ${\Bbb B}_\pi$ is an orthonormal basis for $\pi$, $R(f) \phi$ denotes the action of
$f \in C_c^\infty(\underline{Z}({\Bbb A}_F)\backslash \underline{G}({\Bbb A}_F), \chi)$ on the space of automorphic forms (on
which the center operates via $\chi$) 
by convolution, and $ P(R(f) \phi)$ denotes the period integral of
$R(f) \phi$ along $\underline{H}$:
$$ P(R(f) \phi) = \int_{\underline{H}(F)\underline{Z}({\Bbb A}_F) 
\backslash \underline{H}({\Bbb A}_F)}(R(f) \phi)(h)
\chi(h)^{-1} dh.$$
\end{theorem}
\vspace{.5 cm} 
\noindent{\bf Proof of Theorem 4.1:} The proof of this theorem is via the method of 
relative trace formula, and is a variation on the proof of Hakim-Murnaghan which is
for $\underline{G}= GL_n$, and $\underline{H}$ the subgroup of the fixed points of an involution. We
repeat most of their proof, and observe that the choice of their $\underline{G},\underline{H}$ is 
not really used in the proof, their explicit proof being replaced by 
the general lemmas 4.2 and 4.3 from the theory of Algebraic groups.

It suffices to show  that there exists $f = \otimes f_v \in 
C_c^\infty(\underline{Z}({\Bbb A}_F)\backslash  \underline{G}({\Bbb A}_F), \chi)$  such that $
\Prod_v \Phi_v(f_v,\gamma)$
is nonzero for exactly one $\gamma \in \underline{H}(F)\backslash \underline{G}(F)$. 
We will do this by
choosing $f_v$ appropriately at each place $v$ of $F$. For places in 
$S$, we choose
$f_v$ to be a matrix coefficient of $\pi_v$. 
Since the left and right translates
 of a matrix coefficient 
is  again a matrix coefficient, lemma 4.4 allows us to assume that 
$\Phi_v(f_v, 1) \not = 0$. For places outside $T$, define $f_v(zk)= 
\chi_v(z) $ where $z \in \underline{Z}(F_v)$, and $k$ 
an element of the hyperspecial maximal compact subgroup $\underline{G}(\cO_v)$ of 
$\underline{G}(F_v)$, and define $f_v$ to be 0 outside $\underline{Z}(F_v) \cdot \underline{G}(\cO_v)$; this
is well-defined because $\chi_v$ is unramified. 
For the remaining places of $F$ we choose $f_v \in
C_c^\infty(\underline{Z}(F_v)\backslash \underline{G}(F_v),\chi_v)$ arbitrarily so that 
$\Phi(f_v, 1)$ is nonzero.
In this last choice, ensure that at some place, say an infinite place 
$v_0$ of $F$,
the support of the function  $\Phi(f_{v_0}, g)$     is a certain 
small neighborhood of the identity in $\underline{G}(F_{v_0})/\underline{H}(F_{v_0})$. Since 
the `orbital integral' map is easily
seen to be a surjection from 
$ C_c^\infty( \underline{Z}(F_{v_0})\backslash  \underline{G}(F_{v_0}),\chi_{v_0})$
to  $ C_c^\infty(  \underline{H}(F_{v_0} \backslash  
\underline{G}(F_{v_0}),\chi_{v_0})$, this is possible
to achieve. Finally, we see that $\Prod_v \Phi_v(f_v,\gamma),$
can be assumed to be nonzero for exactly one $\gamma \in \underline{H}(F)\backslash 
\underline{G}(F)$, in fact for $\gamma = 1$, using
lemma 4.2, and the next lemma whose trivial proof we will omit.
This completes the proof of the theorem.

\begin{lemma} Let $X$ and $Y$ be locally compact Hausdorff topological spaces,
and $D \subset X {\times} Y $ be a discrete subset containing a point
$(x_0,y_0)$. Let $K \subset X$ be a compact set containing $x_0$. Then
there exists an open set $U \subset Y$, such that $D \cap (K {\times} U) = 
(x_0,y_0)$.
\end{lemma}

\noindent{\bf Remark :} Because the globalisation theorem is useful in a
variety of contexts, we have tried to state theorem 4.1 quite  generally, 
without restricting $\underline{H}$ to be reductive. In particular this result
contains as a special case the well-known result that locally generic 
supercuspidal representations can be embedded as components of a generic automorphic representation. 
It can also be used to globalise local representations with particular
Fourier-Jacobi, or Bessel models. We will discuss an example involving an
$\underline{H}$ which is neither unipotent, nor reductive.

\vspace{5mm}
\noindent{\bf Example :} Let $G = GL_4 {\times} GL_2$, $H = \Delta GL_2 \cdot N$
where $\Delta GL_2$ is the embedding of $GL_2$ in $G$ given by 
$$g \in GL_2 \rightarrow \left ( \begin{array}{cc} g & 0 \\ 0 & g \end{array}
\right ) {\times} g \in GL_4 {\times} GL_2,$$
and $N$ is the unipotent radical of the standard $(2,2)$ parabolic in $GL_4$,
$$N = \left \{  \left ( \begin{array}{cc} I & X \\ 0 & I \end{array}
\right ) \mid X \in M_2(k)  \right \}.$$
An additive character $\psi_0$ of $k$ defines an additive character $\psi$
of $N$ by $\psi(X) = \psi_0({\rm tr}X).$ This extends to a character $\chi$
of $H$ by $\chi(g,n) = \psi(n).$ A $\chi$-distinguished representation
of $GL_4 {\times} GL_2$ of the form $\pi_1 \otimes \pi_2$ is nothing but a 
representation $\pi_1$ of $GL_4$ whose maximal quotient $\pi_{1,\psi}$ on
which $N$ operates via $\psi$ which is a representation space of $GL_2$ 
has $\pi^\vee_2$ in its quotient. Theorem 4.1 allows us to globalise this to
construct nonvanishing  period integrals of the form: 

$$\int_{GL_2(F){\Bbb A}_F^{\times}\backslash GL_2({\Bbb  A}_F)  
{\times} M_2(F)\backslash M_2({\Bbb A}_F)  }f_1
\left( \begin{array}{cc}  g & gX \\
0 & g \end{array} \right ) f_2(g) \psi(X) dg dX $$
where $f_1$ is a cusp form  belonging to an 
automorphic representation $\Pi_1$ on $GL_4({\Bbb A}_F)$ with $P$ 
as the $(2,2)$ maximal parabolic; 
the function $f_2$ belongs to a cuspidal 
automorphic representation $\Pi_2$ of $GL_2({\Bbb A}_F)$. 
If the non-vanishing of this 
integral can be related to the non-vanishing of the central critical value 
$L(\Lambda^2 \Pi_1 \otimes \Pi_2^{\times}, \frac{1}{2}),$ as is expected, then
following the method of the next section, we will obtain a local 
result about the constituents of the twisted Jacquet
functor $\pi_{1,\psi}$, as conjectured in \cite{deg-Whittaker} in terms of epsilon factors,
and verified there for nonsupercuspidal representations of $GL_4$.

\section{A local consequence}

Let  ${\Bbb K}$ be a commutative semisimple
 cubic algebra over a non-archimedean
local field $k.$  In this section we  study conditions under which
irreducible, admissible representations of $GL_2({\Bbb K})$ have
$GL_2(k)$-invariant linear forms. The results are expressed in terms
 of certain epsilon factors. The case ${\Bbb K} = k 
\oplus k \oplus k$, has been previously studied in \cite{p-trilinear}
 in the odd residue 
characteristic and completed in \cite{prasad_saito}. In fact the method of proof given
in \cite{prasad_saito} works also for the case  when
${\Bbb K}$ is of the form $K \oplus k$ where $K$ is a quadratic field
extension of $k$, but does not seem to work when ${\Bbb K}$
is a cubic field extension of $k$. 
 For the case when ${\Bbb K}= K \oplus k$ with $K$ 
a quadratic field extension of $k$, or ${\Bbb K}$ a cubic field extension 
of $k$,
such a result was proved in \cite{dprasad_amj} exactly for 
those  cases for which the representation 
of   $GL_2({\Bbb K})$ was not supercuspidal;  the local methods employed there 
were quite inadequate to handle supercuspidal representations in these cases.
Our present method proves
such a result in complete generality.

We begin by fixing some notation.
Let $D_k$ be the unique quaternion division algebra over $k$, and let
$D_{{\Bbb K}} = D_k \otimes_k {\Bbb K}.$  If $\pi$  is a discrete series representation
of ${\rm GL}_2({\Bbb K})$ 
(by which we will always mean an irreducible representation which has a
twist whose matrix coefficients are square integrable modulo centre),
we associate an irreducible, admissible representation $\pi'$ of
$D^{\times}_{{\Bbb K}}$ as follows.  If ${\Bbb K} = K \oplus k,$ where $K$  is a
quadratic field extension of $k,$ then ${\rm GL}_2({\Bbb K}) = {\rm GL}_2(K) {\times} {\rm GL}_2(k)$
and the discrete series representation $\pi$ of ${\rm GL}_2({\Bbb K})$  is the tensor
product $\pi_1 \otimes \pi_2$  of a discrete series representation $\pi_1$ of
${\rm GL}_2(K)$ and a discrete series representation $\pi_2$ of ${\rm GL}_2(k).$ In this
case $D^{\times}_{{\Bbb K}} = {\rm GL}_2(K) {\times} D^{\times}_k.$  We define the representation $\pi'$
of $D^{\times}_{{\Bbb K}}$  to be $\pi_1 \otimes \pi'_2$ where $\pi'_2$ is the representation
of $D^{\times}_k$  associated to the discrete series representation $\pi_2$ of
${\rm GL}_2(k)$ by the Jacquet-Langlands correspondence. If ${\Bbb K}$  is a cubic
field extension of $k$ then $D_{{\Bbb K}}$ is the unique quaternion division
algebra over the field ${\Bbb K}$ and we let $\pi'$ be the representation of
$D^{\times}_{{\Bbb K}}$  associated to the discrete series representation $\pi$ of
${\rm GL}_2({\Bbb K})$ (by the Jacquet-Langlands correspondence).

In this paper we will be dealing only with {\em generic}
representations of ${\rm GL}_2({\Bbb K})$ by which we will mean any infinite
dimensional representation of ${\rm GL}_2({\Bbb K})$  if ${\Bbb K}$  is a cubic field
extension  of $k,$  and  if ${\Bbb K} = K \oplus  k$ to be the tensor product of
an infinite dimensional representation of ${\rm GL}_2(K)$  with an infinite
dimensional representation of ${\rm GL}_2(k).$

The main results
proved in \cite{dprasad_amj} were the following.

\begin{enumerate}
\item {\bf Multiplicity one  theorem:} 
For an irreducible, admissible representation
$\pi$ of ${\rm GL}_2({\Bbb K}),$  the space of ${\rm GL}_2(k)$-invariant linear forms on
$\pi$ is at most one-dimensional.

\item{\bf Dichotomy principle:} 
Let  $\pi$  be an irreducible, admissible, generic
representation of ${\rm GL}_2({\Bbb K})$  such that the central character of $\pi$
restricted to $k^{\times} \subseteq  {\Bbb K}^{\times}$  is trivial. Then either there
exists a ${\rm GL}_2(k)$-invariant linear form on $\pi$ which is unique up to
scalars  or  the representation $\pi$ of ${\rm GL}_2({\Bbb K})$  is a discrete
series representation and there exists a $D_k^{\times}$-invariant linear
form on the representation $\pi'$ of $D_{{\Bbb K}}^{\times}$  which is also unique
up to scalars. Moreover, only one of the two possibilities occurs.

\item {\bf Theorem about epsilon factors:} 
Let $\pi$ be an irreducible, admissible, generic
representation of $GL_2({\Bbb K})$ 
such that the central character of $\pi$
restricted to $k^{\times} \subseteq  {\Bbb K}^{\times}$  is trivial. Then for  $\psi$ a
non-trivial additive character of $k$, ~  $\epsilon (M_{{\Bbb K}}^k
\sigma_{\pi},\psi) \cdot \omega_{{\Bbb K} /k} (-1)$ is independent of $\psi$ and
takes the value $+1$ if and only if the representation $\pi$  of ${\rm GL}_2({\Bbb K})$  
has a ${\rm GL}_2(k)$-invariant linear form, and takes the value $-1$ if and only 
if the representation $\pi$ of $GL_2({\Bbb K})$ is a discrete series representation and
the representation $\pi'$ of $D^{\times}_{{\Bbb K}}$ has a $D^{\times}_k$-invariant linear
form.

\end{enumerate}

In the theorem on epsilon factors, 
$M_{{\Bbb K}}^k \sigma_{\pi}$  is a certain 8-dimensional representation
 of the Deligne-Weil group of $k$  defined in \cite{dprasad_amj}, which, 
for ${\Bbb K}$ a cubic field extension of $k$, is
associated  by the process of `tensor induction'
to  the 2-dimensional reprepresentation $\sigma_{\pi}$ of the 
Deligne-Weil group of ${\Bbb K}$ where $\sigma_{\pi}$ is the
Langlands parameter of the  irreducible admissible representation  $\pi$ of
${\rm GL}_2({\Bbb K})$. 
We have used $\omega_{{\Bbb K}/k}$ for the following quadratic
character of $k^{\times}$. 

(a) If ${\Bbb K} = K \oplus k$ with $K$ a quadratic field extension of $k$ then
 $\omega_{{\Bbb K}/k} = \omega_{K/k}$ where $\omega_{K/k}$ is the quadratic character 
 of $k^{\times}$  associated by local
classfield theory to $K$. 

(b) If ${\Bbb K}$ is a 
cubic field extension of $k$ then $\omega_{{\Bbb K}/k}$ will be the trivial 
character of $k^{\times}$ if ${\Bbb K}$ is Galois over $k$, and will be the quadratic character 
$\omega_{L/k}$ of $k^{\times}$ if ${\Bbb K}$ is not Galois over $k$ and $L$ is the unique 
quadratic extension of $k$ contained in the Galois closure of ${\Bbb K}$.

\vspace{2mm}

Of these three theorems, only the first two were proved in complete generality 
in \cite{dprasad_amj}, and
the third theorem on epsilon factors was checked  case-by-case only for non-supercuspidal
representations.

\vspace{3mm}

\noindent{\bf Remark :} In \cite{dprasad_amj}, the epsilon factor used 
was that of an eight dimensional representation 
$M_{{\Bbb K}}^k \sigma_{\pi}$
of the Weil-Deligne
group $W_k'$. In the present work, this Galois theoretic epsilon factor
is changed to the one associated to $\pi$ by the Langlands-Shahidi
method, as in \cite{shahidi_annals127}. 
It is natural to expect that the two epsilon factors are the same, but it has not been
proven yet except in the case $E=F\oplus F\oplus F$, which has been done 
by Ramakrishnan \cite{ramakrishnan_annals152}.
% It would be convenient for us to define the Langlands-Shahidi $\epsilon$-factor 
% as $\epsilon_3(s,\pi)$, so that, one has $$ \epsilon_3(s,\pi) = 
%  \epsilon (M_{{\Bbb K}}^k \sigma_{\pi},\psi) \cdot \omega_{{\Bbb K} /k} (-1).$$

\vspace{.7 cm}

\noindent{\bf Proof of the Theorem on epsilon factors:} The proof of the
theorem on epsilon factors is exactly along the same lines as given in 
\cite{prasad_saito}, except that we globalise a locally distinguished 
representation by means of theorem 4.1, thus via the methods of relative
trace formula instead of the Burger-Sarnak principle that we employed in 
\cite{prasad_saito}.

Given a commutative semisimple cubic
algebra ${\Bbb K}$ over $k$, realise this as the local component of a cubic field
extension $E$ of a number field $F$ at a place $v_0$, i.e., $k = F_{v_0}$, and 
${\Bbb K} = E \otimes_F F_{v_0}$. We assume that both $F$ and $E$ are 
totally real. Assume that $\pi $ has a $GL_2(k)$-invariant 
linear form.  
For $\pi$ non-supercuspidal, the theorem on epsilon factors is proved by
direct calculation in \cite{dprasad_amj}. We therefore assume in the rest of
the proof that $\pi$ is a supercuspidal representation of $GL_2({\Bbb K})$.
Then by theorem 4.1, $\pi$ can be realised as the local component
of an automorphic representation $\Pi$ of $GL_2({\Bbb A}_E)$, and 
there exists an $f \in \Pi$ 
such that 
$$\int_{GL_2(F){\Bbb A}_F^{\times}\backslash GL_2({\Bbb A}_F)} f(h) dh \not = 0,$$
with the further property that $\Pi$ is unramified at any finite place of $E$ 
which is not lying over $v_0$. By theorem 1.1, this means that $L_3(\frac{1}{2},\Pi)
\not = 0.$ Therefore, in particular, the sign 
$\epsilon_3(\frac{1}{2},\Pi)$ in the functional equation is 1.

Before proceeding further, define a character $\omega_{E/F}: 
{\Bbb A}_F^{\times} \rightarrow
{\pm 1}$ to be the sign of the permutation representation of the Galois group
Gal$(\bar{F}/F)$ on the set Gal$(\bar{F}/F) / {\rm Gal}(\bar{F}/E)$. Clearly
$\omega_{E/F}$ has local components $\omega_{{\Bbb {\Bbb K}} /k}$ 
where $k$ is a completion
of $F$, and ${\Bbb K} = k \otimes_F E$. In particular, 
$\Prod_v \omega_{E_v /F_v}(-1) =1$
where the product is over all the places $v$ of $F$, and 
$E_v = E \otimes_F F_v$.

Since $$\epsilon_3( \frac{1}{2}, \Pi) 
= \prod_v \epsilon_3( \frac{1}{2}, \Pi_v ) =\prod_v \epsilon_3( \frac{1}{2}, \Pi_v )
\omega_{E_v /F_v}(-1) =1   ,$$
(where $v$ belongs to the set of places of $F$), 
if we can prove that all the  factors $\epsilon_3( \frac{1}{2}, \Pi_v )
\omega_{E_v /F_v}(-1) $ 
are equal to  1, except at $v_0$, we would be able to deduce that
$\epsilon_3( \frac{1}{2}, \Pi_{v_0} )  \omega_{{\Bbb K} / k}(-1) = 1$ 
too. But
$\epsilon_3( \frac{1}{2}, \Pi_v ) \omega_{E_v /F_v}(-1) = 1$ 
for all finite places $v$, $v \not = v_0$, 
as $\Pi_v$ is an unramified representation,
for which this is an easy calculation as
done in  proposition 8.7 of \cite{dprasad_amj}.  (As an erratum to 
\cite{dprasad_amj},
we note that the factor 
$\omega_{E_v /F_v}(-1)$ is missing in the statement of the proposition 
8.7 of \cite{dprasad_amj}.)
For $v$ Archimedean, again 
$\epsilon_3( \frac{1}{2}, \Pi_v )  
\omega_{E_v /F_v}(-1)= 1$. (We note that since we 
have assumed that all the real places of $F$ split in $E$, $\omega_{E_v /F_v}(-1)= 1$, and 
that by the theorem on triple products at infinity, cf. 
\cite{p-trilinear}, 
$\epsilon_3( \frac{1}{2}, \Pi_v )  = 1$.)
Thus we conclude 
that if $\pi_{v_0}$ has a $GL_2(k)$-invariant linear form, then
$\epsilon_3( \frac{1}{2}, \Pi_{v_0} )  \omega_{{\Bbb K} / k}(-1) = 1$. 

On the
other hand if $\pi$ does not have a $GL_2(k)$-invariant linear form,
then by the dichotomy principle proved in \cite{dprasad_amj}, $\pi$ 
is a discrete series representation, and  there exists a $D_k^{\times}$-invariant 
linear
form on the representation $\pi'$ of $D_{{\Bbb K}}^{\times}$ associated to $\pi$
by the Jacquet-Langlands correspondence. Globalise $\pi'$ to a distinguished
automorphic form $\Pi'$ on $(D\otimes_F E)^{\times}$ where  
$D$ is the quaternion division 
algebra over $F$, which remains a division algebra at the place $v_0$,  and 
which is a division algebra at exactly one place at infinity 
of $F$, but at no other finite or infinite place of $F$ besides these two.
We can assume again that $\Pi'$ is unramified at all finite places of $E$ 
except the one over $v_0$. This time $\epsilon_3( \frac{1}{2}, \Pi_v )  
\omega_{E_v /F_v}(-1)= -1$ at the infinite place of $F$ where
$D$ has remained a division algebra, and is equal to 1
 at all the other
places besides $v_0$. Thus, $\epsilon_3( \frac{1}{2}, \Pi_{v_0} )  
\omega_{{\Bbb K} / k }(-1)= -1$.
This completes the proof of the theorem on epsilon factors.

\vspace{0.2cm}

\noindent Dipendra Prasad

\noindent School of Mathematics, Tata Institute of Fundamental Research,

\noindent Colaba, Mumbai-400005, INDIA

\noindent Email: dprasad@math.tifr.res.in

\vspace{0.2cm}
Rainer Schulze-Pillot\\
Fachrichtung 6.1 Mathematik,
Universit\"at des Saarlandes (Geb. E2.4)\\
Postfach 151150\\
66041 Saarbr\"ucken, Germany\\
email: schulzep@math.uni-sb.de

\end{document}